\documentclass[12pt]{amsart}
\usepackage[utf8]{inputenc}
\usepackage{amsmath,amssymb,amsthm}
\usepackage{enumerate, euler, palatino}
\usepackage[colorlinks=true]{hyperref}
\hypersetup{urlcolor=blue, citecolor=red , linkcolor= blue}
\usepackage{cite}
\usepackage[a4paper,margin=2.5cm]{geometry}

\theoremstyle{plain}
\newtheorem{theorem}{Theorem}[section]
\newtheorem{proposition}[theorem]{Proposition}
\newtheorem{lemma}[theorem]{Lemma}
\newtheorem{corollary}[theorem]{Corollary}
\theoremstyle{definition}

\newtheorem{remark}[theorem]{Remark}

\title[Singular values of iterated Aluthge transforms]{Singular values of iterated Aluthge transforms of compact operators}
\author{Hranislav Stankovi\'c}
\address{Faculty of Electronic Engineering, University of Ni\v{s}, Aleksandra Medvedeva 4, 18000 Ni\v{s}, Serbia}
\email{hranislav.stankovic@elfak.ni.ac.rs}

\begin{document}
	
	\begin{abstract}
		Yamazaki proved that the norms of the iterated Aluthge transforms of an operator converge to its spectral radius, and Ando showed that for matrices the singular values of the iterates converge to the moduli of the eigenvalues. We prove that the latter holds for every compact operator on a Hilbert space and every generalized Aluthge transform $\Delta_t$ with $0<t<1$: for every $j$, the $j$-th singular value of the $n$-th iterate converges to $|\lambda_j(T)|$, and the products of the $k$ largest singular values, as well as the Ky Fan norms, decrease monotonically to their limits. This is the Aluthge counterpart of the Yamamoto--Davis theorem. We also determine the limit in the Schatten $p$-norm whenever $T$ belongs to a strictly smaller Schatten class.
	\end{abstract}
	
	\maketitle
	\bigskip
	\section{Introduction}
	\bigskip
	
	Let $\mathcal{H}$ be a complex Hilbert space and let $\mathfrak{B}(\mathcal{H})$ be the algebra of bounded linear operators on $\mathcal{H}$, with the operator norm $\|\cdot\|$ and the identity $I$. For $T\in\mathfrak{B}(\mathcal{H})$ we write $T^*$, $\ker T$, $\operatorname{ran}T$, $\sigma(T)$ and $r(T)$ for the adjoint, the kernel, the range, the spectrum and the spectral radius of $T$, and $|T|=(T^*T)^{1/2}$ for its modulus. The operator $T$ is positive if $\langle Tx,x\rangle\ge0$ for every $x\in\mathcal{H}$, and normal if $T^*T=TT^*$. For positive $T$ the powers $T^{\alpha}$, $\alpha>0$, are defined by the continuous functional calculus. Every $T\in\mathfrak{B}(\mathcal{H})$ has a unique polar decomposition $T=U|T|$ in which $U$ is a partial isometry with initial space $(\ker T)^\perp$. For $t\in[0,1]$ the generalized Aluthge transform of $T$ is
	\begin{equation*}
		\Delta_t(T)=|T|^{t}U|T|^{1-t},
	\end{equation*}
	where $|T|^{0}=I$. For $t=1/2$ this is the Aluthge transform, introduced by Aluthge \cite{Aluthge} and studied systematically by Jung, Ko and Pearcy \cite{JKP}; the general case is due to Okubo \cite{Okubo} (see also \cite{AMS}). At the endpoints, $\Delta_0(T)=T$, while $\Delta_1(T)=|T|U$ is the Duggal transform of $T$ \cite{FJKP}. In the literature $\Delta_t$ is usually called the $\lambda$-Aluthge transform, and we denote its parameter by $t$ only in order to keep $\lambda$ for eigenvalues. The iterates are defined by $\Delta_t^0(T)=T$ and $\Delta_t^{n+1}(T)=\Delta_t(\Delta_t^n(T))$. Two facts are basic. The first is that $\Delta_t(T)$ has the same spectrum as $T$. For $t=1/2$ this is due to Jung, Ko and Pearcy \cite{JKP}, and their argument works for every $t\in[0,1]$: the factorizations $T=\big(U|T|^{1-t}\big)|T|^{t}$ and $\Delta_t(T)=|T|^{t}\big(U|T|^{1-t}\big)$ give $\sigma(T)\setminus\{0\}=\sigma(\Delta_t(T))\setminus\{0\}$, and $T$ is invertible if and only if $\Delta_t(T)$ is. The second is that $\Delta_t$ does not increase the operator norm,
	\begin{equation}\label{eq:contraction}
		\|\Delta_t(T)\|\le\|T\|,\qquad T\in \mathfrak{B}(\mathcal{H}),\quad 0\le t\le1,
	\end{equation}
	which is immediate from $\|U\|\le1$, $\||T|^{t}\|=\|T\|^{t}$ and $\||T|^{1-t}\|=\|T\|^{1-t}$.
	
	Yamazaki \cite{Yamazaki} showed that for $t=1/2$ the decreasing sequence of norms of the iterates ends exactly at the spectral radius. The same is true for every $t\in(0,1)$: this was proved for matrices in \cite{AMS} and for invertible operators by Tam \cite{Tam}, and in general it is the one-variable case of a theorem of Benhida, Curto, Lee and Yoon \cite{BCLY} on the generalized spherical Aluthge transform of a commuting tuple. Thus
	\begin{equation}\label{eq:yamazaki}
		\lim_{n\to\infty}\|\Delta_t^n(T)\|=r(T),\qquad T\in \mathfrak{B}(\mathcal{H}),\quad 0<t<1.
	\end{equation}
	This is an analogue of Gelfand's formula
	\begin{equation}\label{eq:gelfand}
		r(T)=\lim_{n\to\infty}\|T^n\|^{1/n},\qquad T\in \mathfrak{B}(\mathcal{H}),
	\end{equation}
	with powers replaced by Aluthge iterates. The restriction to $0<t<1$ is essential. For $t=0$ the iterates are constant, and for the Duggal transform $\Delta_1$ the formula fails, as an example in \cite{BCLY} shows. The Duggal iterates of an arbitrary operator are described in Proposition \ref{prop:duggal} below. See also \cite{Wang} for a simplified proof of the case $t=1/2$.
	
	Gelfand's formula has a well known refinement due to Yamamoto \cite{Yamamoto}. If $T$ is a $d\times d$ matrix with singular values $s_1(T)\ge\cdots\ge s_d(T)$ and eigenvalues $\lambda_1(T),\ldots,\lambda_d(T)$ arranged so that $|\lambda_1(T)|\ge\cdots\ge|\lambda_d(T)|$, then
	\begin{equation}\label{eq:yamamoto}
		\lim_{n\to\infty}s_j(T^n)^{1/n}=|\lambda_j(T)|,\qquad j=1,\ldots,d.
	\end{equation}
	For $j=1$ this is \eqref{eq:gelfand}. Yamamoto \cite{Yamamoto2} obtained infinite-dimensional versions of \eqref{eq:yamamoto} in special cases, and Davis \cite{Davis} proved that \eqref{eq:yamamoto} holds verbatim for every compact operator on a Hilbert space, a result now known as the Yamamoto--Davis theorem; a recent spatial strengthening is given in \cite{NS}. The purpose of this paper is to prove the corresponding refinement of \eqref{eq:yamazaki}. For a compact operator $T$ we write $s_j(T)$ for the singular values and $\lambda_j(T)$ for the eigenvalues, with the conventions recalled in Section \ref{sec:prelim}.
	
	\begin{theorem}\label{thm:main-intro}
		Let $T$ be a compact operator on a complex separable Hilbert space $\mathcal{H}$ and let $0<t<1$. Then
		\begin{equation*}
			\lim_{n\to\infty}s_j(\Delta_t^n(T))=|\lambda_j(T)|
		\end{equation*}
		for every $j$. Moreover, for every $k$ the products $\prod_{j=1}^{k}s_j(\Delta_t^n(T))$ decrease monotonically to $\prod_{j=1}^{k}|\lambda_j(T)|$ as $n\to\infty$.
	\end{theorem}
	
	Theorem \ref{thm:main-intro} is the Aluthge counterpart of \eqref{eq:yamamoto}: the powers $T^n$ are replaced by the iterates $\Delta_t^n(T)$, no normalization by the exponent is needed, and the convergence is monotone at every level $k$. It also bears on a question about the iterates themselves. These need not converge, not even in the weak operator topology \cite{CJL}, and in \cite[Section 3.3]{APS3} Antezana, Pujals and Stojanoff ask for a characterization of the operators whose iterates converge, suggesting compact operators as the first case to study. Theorem \ref{thm:main-intro} shows that for compact operators the singular values converge whether or not the iterates do.
	
	The proof of \eqref{eq:yamamoto} rests on two facts: the antisymmetric tensor power $\wedge^k$ turns products of the $k$ largest singular values, respectively eigenvalue moduli, into the norm, respectively the spectral radius, of $\wedge^kT$, and $\wedge^k$ commutes with powers. For matrices the Aluthge analogue of this argument was carried out by Ando \cite{Ando}: compound matrices commute with the Aluthge transform, and \eqref{eq:yamazaki} applied to them gives Theorem \ref{thm:main-intro} for $t=1/2$ and $d<\infty$. The same argument with the matrix case of \eqref{eq:yamazaki} gives the matrix case for every $t\in(0,1)$, see \cite[Theorem 1.4]{HT}. Our proof of Theorem \ref{thm:main-intro} follows the same route, with the commutation with powers replaced by the identity
	\begin{equation}\label{eq:commute}
		\wedge^k\Delta_t(T)=\Delta_t(\wedge^kT),
	\end{equation}
	which follows from the compatibility of polar decompositions with antisymmetric tensor powers. Identity \eqref{eq:commute}, combined with \eqref{eq:contraction}, is the mechanism behind the log-majorization
	\begin{equation}\label{eq:logmaj}
		\prod_{j=1}^{k}s_j(\Delta_t(T))\le\prod_{j=1}^{k}s_j(T),\qquad k\in\mathbb{N},
	\end{equation}
	of Antezana, Massey and Stojanoff \cite[Remark 3.12]{AMS}, valid for every compact operator $T$, and iterating \eqref{eq:logmaj} already gives the monotonicity in Theorem \ref{thm:main-intro}. What is new here is the limit for compact operators on infinite-dimensional spaces. Formula \eqref{eq:yamazaki}, applied to $\wedge^kT$ instead of $T$, identifies it as $r(\wedge^kT)$, and the individual singular values are then recovered from ratios of consecutive products. Two ingredients are needed beyond the matrix case: \eqref{eq:yamazaki} in its full generality \cite{BCLY}, and the identification of $r(\wedge^kT)$ with $\prod_{j\le k}|\lambda_j(T)|$ for a compact operator $T$, which for matrices is immediate from the Jordan form but in general requires an argument.

	\bigskip 
	
	The paper is organized as follows. Section \ref{sec:prelim} collects the facts about antisymmetric tensor powers, Section \ref{sec:commute} relates them to the Aluthge transforms, Section \ref{sec:main} contains the proof of the main theorem, and Section \ref{sec:norms} its consequences for Ky Fan norms, for arbitrary unitarily invariant norms in the finite-dimensional case, and for the norm ideals of compact operators and the Schatten norms in general. Section \ref{sec:duggal} describes the iterated Duggal transform, which shows what happens at the excluded endpoint $t=1$.
	
	\bigskip
	\section{Preliminaries}\label{sec:prelim}
	\bigskip
	
	Throughout, $\mathcal{H}$ is a complex separable Hilbert space of dimension $d\in\mathbb{N}\cup\{\infty\}$, and $\mathcal{K}(\mathcal{H})$ denotes the ideal of compact operators. The separability assumption is discussed in Remark \ref{rem:separable} at the end of the paper. By the Riesz theory of compact operators \cite[Theorem VII.7.1]{Conway}, a compact operator has at most countably many non-zero eigenvalues, each of finite geometric multiplicity, and they accumulate at most at $0$. Let $T\in\mathcal{K}(\mathcal{H})$. The singular values of $T$ are the eigenvalues of $|T|$, and the non-zero ones are arranged in the non-increasing sequence $s_1(T)\ge s_2(T)\ge\cdots$, each eigenvalue $s$ of $|T|$ being repeated $\dim\ker(|T|-s)$ times. Similarly, $\lambda_1(T),\lambda_2(T),\ldots$ are the non-zero eigenvalues of $T$ arranged so that $|\lambda_1(T)|\ge|\lambda_2(T)|\ge\cdots$, each repeated according to its algebraic multiplicity. The algebraic multiplicity of a non-zero eigenvalue $\mu$ of $T$ is the dimension of the range of the Riesz idempotent $E_\mu$ of $T$ associated with $\{\mu\}$, which is finite \cite[Corollary VII.7.8]{Conway}, and it equals $\dim\ker(T-\mu)^{m}$ for every sufficiently large $m$ \cite[Chapter VII, Sections 3 and 4]{DS}. For the self-adjoint operator $|T|$ it reduces to $\dim\ker(|T|-\mu)$, so that the two sequences are formed by the same rule. For a compact operator $X$ let $\rho(X)\in\{0,1,2,\ldots\}\cup\{\infty\}$ be the number of non-zero eigenvalues of $X$ counted with algebraic multiplicity. So far the sequence $(\lambda_j(T))_j$ has the indices $j\le\rho(T)$, and $(s_j(T))_j$ the indices $j\le\rho(|T|)$. We extend both by zeros, putting $\lambda_j(T)=0$ for $j>\rho(T)$ and $s_j(T)=0$ for $j>\rho(|T|)$, so that both sequences are indexed by $j\in\{1,\ldots,d\}$ if $d<\infty$ and by $j\in\mathbb{N}$ if $d=\infty$, and $\lambda_j(T)\neq0$ exactly for $j\le\rho(T)$. In finite dimensions $(s_j(T))_j$ then consists of all eigenvalues of $|T|$, and $(\lambda_j(T))_j$ of all eigenvalues of $T$ counted with algebraic multiplicity. In infinite dimensions the eigenvalue $0$ is left out of the count because it need not be an isolated point of $\sigma(T)$, and need not be an eigenvalue at all.
	
	\bigskip 
	
	Let $k\in\mathbb{N}$ with $k\le d$. For $x_1,\ldots,x_k\in\mathcal{H}$ we put
	\begin{equation*}
		x_1\wedge\cdots\wedge x_k=\frac{1}{\sqrt{k!}}\sum_{\pi}\operatorname{sgn}(\pi)\,x_{\pi(1)}\otimes\cdots\otimes x_{\pi(k)},
	\end{equation*}
	the sum being taken over all permutations $\pi$ of $\{1,\ldots,k\}$, so that
	\begin{equation*}
		\langle x_1\wedge\cdots\wedge x_k,\,y_1\wedge\cdots\wedge y_k\rangle=\det\big(\langle x_i,y_j\rangle\big)_{i,j=1}^{k}.
	\end{equation*}
	The $k$-th antisymmetric tensor power $\wedge^k\mathcal{H}$ is the closed linear span of these vectors in the Hilbert space tensor power $\mathcal{H}^{\otimes k}$. If $(e_i)$ is an orthonormal basis of $\mathcal{H}$, then the vectors $e_{i_1}\wedge\cdots\wedge e_{i_k}$ with $i_1<\cdots<i_k$ form an orthonormal basis of $\wedge^k\mathcal{H}$, see \cite[Formula (1.10)]{Simon}. For a closed subspace $\mathcal{L}$ of $\mathcal{H}$ we write $\wedge^k\mathcal{L}$ for the closed linear span of all $x_1\wedge\cdots\wedge x_k$ with $x_1,\ldots,x_k\in\mathcal{L}$. For bounded operators $S_1,\ldots,S_k$ on Hilbert spaces $\mathcal{H}_1,\ldots,\mathcal{H}_k$ the tensor product $S_1\otimes\cdots\otimes S_k$ is the bounded operator on $\mathcal{H}_1\otimes\cdots\otimes\mathcal{H}_k$ determined by $(S_1\otimes\cdots\otimes S_k)(x_1\otimes\cdots\otimes x_k)=S_1x_1\otimes\cdots\otimes S_kx_k$, and it satisfies $\|S_1\otimes\cdots\otimes S_k\|=\|S_1\|\cdots\|S_k\|$, see \cite[Section 2.6]{KR}. For $T\in\mathfrak{B}(\mathcal{H})$ we write $T^{\otimes k}$ for the $k$-fold tensor product of $T$ with itself. The subspace $\wedge^k\mathcal{H}$ is reducing for $T^{\otimes k}$. Indeed, if $W_\pi$ denotes the unitary operator on $\mathcal{H}^{\otimes k}$ which permutes the factors of elementary tensors according to the permutation $\pi$, then $T^{\otimes k}$ commutes with every $W_\pi$, hence with the antisymmetrizer $\frac{1}{k!}\sum_{\pi}\operatorname{sgn}(\pi)W_\pi$, and the latter is the orthogonal projection onto $\wedge^k\mathcal{H}$. We write $\wedge^kT$ for the restriction of $T^{\otimes k}$ to $\wedge^k\mathcal{H}$, so that
	\begin{equation*}
		\wedge^kT\,(x_1\wedge\cdots\wedge x_k)=Tx_1\wedge\cdots\wedge Tx_k.
	\end{equation*}
	It is convenient to set $\wedge^0\mathcal{H}=\mathbb{C}$ and $\wedge^0T=1$. The facts collected in Lemmas \ref{lem:wedge} and \ref{lem:wedge-compact} below are standard, see \cite[Chapter I.5]{Bhatia} for the finite-dimensional case and \cite[Section 1.5]{Simon} in general. We give short arguments or precise references for all of them, except for the last formula in Lemma \ref{lem:wedge-compact}, for which we could not locate a direct proof in the literature and therefore include one.
	
	\begin{lemma}\label{lem:wedge}
		Let $k\in\mathbb{N}$ with $k\le d$.
		\begin{enumerate}[\rm($i$)]
			\item For $T,S\in\mathfrak{B}(\mathcal{H})$ one has
			\begin{equation*}
				\wedge^k(TS)=(\wedge^kT)(\wedge^kS),\qquad(\wedge^kT)^*=\wedge^k(T^*),\qquad\wedge^kI=I.
			\end{equation*}
			\item For $T\in\mathfrak{B}(\mathcal{H})$ one has $\|\wedge^kT\|\le\|T\|^k$.
			\item If $T\in\mathfrak{B}(\mathcal{H})$ is positive, then $\wedge^kT$ is positive and $\wedge^k(T^{1/2})=(\wedge^kT)^{1/2}$.
		\end{enumerate}
	\end{lemma}
	
	\begin{proof}
		The identities $(TS)^{\otimes k}=T^{\otimes k}S^{\otimes k}$, $(T^{\otimes k})^*=(T^*)^{\otimes k}$ and $I^{\otimes k}=I$ hold on all of $\mathcal{H}^{\otimes k}$. Since $\wedge^k\mathcal{H}$ is reducing for each of the operators involved, they pass to the restrictions, which gives ($i$). The estimate in ($ii$) follows from $\|\wedge^kT\|\le\|T^{\otimes k}\|=\|T\|^k$, again because $\wedge^kT$ is a restriction of $T^{\otimes k}$.
		
		Let $T$ be positive. Writing $T=T^{1/2}T^{1/2}$ and using ($i$) we obtain
		\begin{equation*}
			\wedge^kT=\big(\wedge^k(T^{1/2})\big)^*\big(\wedge^k(T^{1/2})\big),
		\end{equation*}
		so that $\wedge^kT$ is positive. Applying this to $T^{1/2}$ in place of $T$ shows that $\wedge^k(T^{1/2})$ is positive as well, and by ($i$) its square is $\wedge^kT$. The positive square root being unique, $\wedge^k(T^{1/2})=(\wedge^kT)^{1/2}$, which proves ($iii$).
	\end{proof}
	
	The square root in Lemma \ref{lem:wedge}($iii$) may be replaced by an arbitrary power. This is what makes the generalized Aluthge transforms accessible.
	
	\begin{lemma}\label{lem:power}
		Let $k\in\mathbb{N}$ with $k\le d$, let $T\in \mathfrak{B}(\mathcal{H})$ be positive and let $t\in(0,1]$. Then $\wedge^k(T^{t})=(\wedge^kT)^{t}$.
	\end{lemma}
	
	\begin{proof}
		The case $t=1$ is trivial and the case $t=1/2$ is contained in Lemma \ref{lem:wedge}($iii$). Iterating the latter gives $\wedge^k(T^{2^{-m}})=(\wedge^kT)^{2^{-m}}$ for every integer $m\ge0$, since the identity for $m$ yields
		\begin{align*}
			\wedge^k\big(T^{2^{-(m+1)}}\big)&=\wedge^k\big((T^{2^{-m}})^{1/2}\big)=\big(\wedge^k(T^{2^{-m}})\big)^{1/2}\\&=\big((\wedge^kT)^{2^{-m}}\big)^{1/2}=(\wedge^kT)^{2^{-(m+1)}}.
		\end{align*}
		A dyadic rational in $(0,1]$ has the form $t=i2^{-m}$ with integers $m\ge0$ and $1\le i\le 2^{m}$. Because $\wedge^k$ is multiplicative and $T^{i2^{-m}}=(T^{2^{-m}})^{i}$, we obtain
		\begin{equation*}
			\wedge^k(T^{t})=\big(\wedge^k(T^{2^{-m}})\big)^{i}=\big((\wedge^kT)^{2^{-m}}\big)^{i}=(\wedge^kT)^{t},
		\end{equation*}
		so the identity $\wedge^k(T^{t})=(\wedge^kT)^{t}$ holds for every dyadic rational $t\in(0,1]$.
		
		Let now $t\in(0,1)$ be arbitrary and choose dyadic rationals $t_1\le t_2\le\cdots$ in $(0,1]$ with $t_j\to t$. Put $M=\max\{1,\|T\|,\|\wedge^kT\|\}$. On $[0,1]$ the functions $u\mapsto u^{t_j}$ decrease pointwise to $u\mapsto u^{t}$, and on $[1,M]$ they increase pointwise to it. In both cases the limit function is continuous, because $t>0$, so Dini's theorem applies on each of the two intervals and yields uniform convergence on $[0,M]$. The spectra of $T$ and of $\wedge^kT$ are contained in $[0,M]$, so the continuous functional calculus turns this uniform convergence into convergence in norm,
		\begin{equation*}
			T^{t_j}\longrightarrow T^{t}\qquad\text{and}\qquad(\wedge^kT)^{t_j}\longrightarrow(\wedge^kT)^{t}.
		\end{equation*}
		The map $S\mapsto\wedge^kS$ is norm continuous on $\mathfrak{B}(\mathcal{H})$, so letting $j\to\infty$ in the identity $\wedge^k(T^{t_j})=(\wedge^kT)^{t_j}$, which holds by the first part of the proof, we obtain $\wedge^k(T^{t})=(\wedge^kT)^{t}$.
	\end{proof}
	
	The next lemma is the link between $\wedge^k$ and the quantities in Theorem \ref{thm:main-intro}.
	
	\begin{lemma}\label{lem:wedge-compact}
		Let $k\in\mathbb{N}$ with $k\le d$ and let $T\in\mathcal{K}(\mathcal{H})$. Then $\wedge^kT$ is compact and
		\begin{equation*}
			\|\wedge^kT\|=\prod_{j=1}^{k}s_j(T),\qquad r(\wedge^kT)=\prod_{j=1}^{k}|\lambda_j(T)|.
		\end{equation*}
	\end{lemma}
	
	\begin{proof}
		The tensor product of two compact operators is compact \cite{ZK}, so an induction on $k$ shows that $T^{\otimes k}$ is compact: the case $k=1$ is the hypothesis, and if $T^{\otimes(k-1)}$ is compact, then so is $T^{\otimes k}=T^{\otimes(k-1)}\otimes T$ under the canonical identification of $\mathcal{H}^{\otimes k}$ with $\mathcal{H}^{\otimes(k-1)}\otimes\mathcal{H}$. Hence $\wedge^kT$ is compact, being the restriction of $T^{\otimes k}$ to a closed subspace. The first formula is proved in \cite[Formula (1.14)]{Simon}.
		
		We turn to the second formula and abbreviate $\gamma=\prod_{j\le k}|\lambda_j(T)|$. Fix $\varepsilon\in(0,1]$ and let $\Lambda_\varepsilon=\{\lambda\in\sigma(T):|\lambda|\ge\varepsilon\}$. By \cite[Theorem VII.7.1]{Conway} this is a finite set of isolated points of $\sigma(T)$, hence open and closed in $\sigma(T)$. Let $E$ be the Riesz idempotent of $T$ associated with $\Lambda_\varepsilon$ and put $\mathcal{N}=\operatorname{ran}E$ and $\mathcal{M}=\ker E$. Here $E=0$ if $\Lambda_\varepsilon$ is empty and $E=I$ if $\Lambda_\varepsilon=\sigma(T)$, and the spectrum of the operator on the space $\{0\}$ is understood to be empty. Then $E$ commutes with $T$, so both subspaces are invariant under $T$, and the operators $A=T|_{\mathcal{N}}$ and $B=T|_{\mathcal{M}}$ satisfy $\sigma(A)=\Lambda_\varepsilon$ and $\sigma(B)=\sigma(T)\setminus\Lambda_\varepsilon$. Moreover, $E$ is the sum of the Riesz idempotents $E_\mu$, $\mu\in\Lambda_\varepsilon$, so that $\mathcal{N}$ is the direct sum of the finite-dimensional spaces $\operatorname{ran}E_\mu$. Hence $\mathcal{N}$ is finite-dimensional, and the eigenvalues of $A$, repeated according to algebraic multiplicity, are precisely those $\lambda_j(T)$ with $|\lambda_j(T)|\ge\varepsilon$. For these facts we refer to \cite[Section VII.4 and Corollary VII.7.8]{Conway}. In particular $r(B)\le\varepsilon$, and $\|A\|\le\|T\|$ because $A$ is a restriction of $T$.
		
		The subspaces $\mathcal{N}$ and $\mathcal{M}$ need not be orthogonal in $\mathcal{H}$, so we pass to the external orthogonal direct sum $\mathcal{N}\oplus\mathcal{M}$, the Hilbert space of all pairs $(y,z)$ with $y\in\mathcal{N}$, $z\in\mathcal{M}$ and $\|(y,z)\|^2=\|y\|^2+\|z\|^2$. The map $V\colon\mathcal{H}\to\mathcal{N}\oplus\mathcal{M}$ given by $Vx=(Ex,(I-E)x)$ is bounded and bijective, its inverse $(y,z)\mapsto y+z$ is bounded as well, and $VTV^{-1}=A\oplus B$. The restriction of $V^{\otimes k}$ to $\wedge^k\mathcal{H}$ is therefore a bounded invertible map onto $\wedge^k(\mathcal{N}\oplus\mathcal{M})$ which intertwines $\wedge^kT$ and $\wedge^k(A\oplus B)$, whence
		\begin{equation*}
			r(\wedge^kT)=r\big(\wedge^k(A\oplus B)\big).
		\end{equation*}
		
		Choose orthonormal bases $(f_i)$ of $\mathcal{N}$ and $(g_i)$ of $\mathcal{M}$. Together they form an orthonormal basis of $\mathcal{N}\oplus\mathcal{M}$, and the wedge products of $k$ distinct members of it form an orthonormal basis of $\wedge^k(\mathcal{N}\oplus\mathcal{M})$. For $\ell\in\{0,1,\ldots,k\}$ let $\mathcal{D}_\ell$ be the closed span of those basis vectors having exactly $\ell$ factors from $(f_i)$, so that $\mathcal{D}_\ell=\{0\}$ when $\ell>\dim\mathcal{N}$ or $k-\ell>\dim\mathcal{M}$. Since $A\oplus B$ leaves $\mathcal{N}$ and $\mathcal{M}$ invariant, each $\mathcal{D}_\ell$ is invariant under $\wedge^k(A\oplus B)$, and the identification
		\begin{equation*}
			f_{i_1}\wedge\cdots\wedge f_{i_\ell}\wedge g_{j_1}\wedge\cdots\wedge g_{j_{k-\ell}}\longmapsto(f_{i_1}\wedge\cdots\wedge f_{i_\ell})\otimes(g_{j_1}\wedge\cdots\wedge g_{j_{k-\ell}})
		\end{equation*}
		carries the restriction of $\wedge^k(A\oplus B)$ to $\mathcal{D}_\ell$ onto $\wedge^{\ell}A\otimes\wedge^{k-\ell}B$, where $\wedge^0$ acts as the identity on $\wedge^0\mathcal{N}=\mathbb{C}$ and on $\wedge^0\mathcal{M}=\mathbb{C}$. Consequently $\wedge^k(A\oplus B)$ is the orthogonal direct sum of the operators $\wedge^{\ell}A\otimes\wedge^{k-\ell}B$ over those $\ell$ for which $\mathcal{D}_\ell\neq\{0\}$, and since the spectrum of a finite orthogonal direct sum is the union of the spectra,
		\begin{equation*}
			r(\wedge^kT)=\max_{\mathcal{D}_\ell\neq\{0\}}r\big(\wedge^{\ell}A\otimes\wedge^{k-\ell}B\big).
		\end{equation*}
		
		Since $\|S\otimes R\|=\|S\|\,\|R\|$ for bounded operators $S$ and $R$, we have $\|(S\otimes R)^n\|=\|S^n\otimes R^n\|=\|S^n\|\,\|R^n\|$, and Gelfand's formula gives $r(S\otimes R)=r(S)r(R)$. If $\ell<k$, then $(\wedge^{k-\ell}B)^n=\wedge^{k-\ell}(B^n)$ and Lemma \ref{lem:wedge}($ii$), applied on $\mathcal{M}$, give
		\begin{equation*}
			r(\wedge^{k-\ell}B)=\lim_{n\to\infty}\big\|\wedge^{k-\ell}(B^n)\big\|^{1/n}\le\lim_{n\to\infty}\|B^n\|^{(k-\ell)/n}=r(B)^{k-\ell}\le\varepsilon^{k-\ell},
		\end{equation*}
		while $r(\wedge^{\ell}A)\le\|\wedge^{\ell}A\|\le\|T\|^{\ell}$. Since $\varepsilon\le1$, every summand with $\ell<k$ therefore has spectral radius at most $\|T\|^{\ell}\varepsilon^{k-\ell}\le c\varepsilon$, where $c=\max\{1,\|T\|^{k}\}$.
		
		Suppose first that $\gamma=0$, that is, $\rho(T)<k$. Then $\dim\mathcal{N}\le\rho(T)<k$ for every $\varepsilon\in(0,1]$, so the summand with $\ell=k$ is absent and $r(\wedge^kT)\le c\varepsilon$. Letting $\varepsilon\to0$ gives $r(\wedge^kT)=0=\gamma$. Suppose now that $\gamma>0$ and choose $\varepsilon<|\lambda_k(T)|$. Then $\dim\mathcal{N}\ge k$ and the $k$ largest eigenvalue moduli of $A$ are $|\lambda_1(T)|,\ldots,|\lambda_k(T)|$. The eigenvalues of $\wedge^kA$ are the products $\lambda_{i_1}(A)\cdots\lambda_{i_k}(A)$ with $i_1<\cdots<i_k$ \cite[Chapter I.5]{Bhatia}, and the largest of their moduli is $\gamma$, so the summand with $\ell=k$ has spectral radius exactly $\gamma$. Hence
		\begin{equation*}
			\gamma\le r(\wedge^kT)\le\max\{\gamma,c\varepsilon\},
		\end{equation*}
		and letting $\varepsilon\to0$ gives $r(\wedge^kT)=\gamma$.
	\end{proof}
	
	The second formula of Lemma \ref{lem:wedge-compact} is stated without proof in \cite{Davis}, where it serves exactly the same purpose as here. It is closely tied to the Yamamoto--Davis theorem. Since $(\wedge^kT)^{n}=\wedge^k(T^{n})$ by Lemma \ref{lem:wedge}($i$), the first formula of Lemma \ref{lem:wedge-compact} and Gelfand's formula give
	\begin{equation*}
		r(\wedge^kT)=\lim_{n\to\infty}\big\|\wedge^k(T^{n})\big\|^{1/n}=\lim_{n\to\infty}\Big(\prod_{j=1}^{k}s_j(T^{n})\Big)^{1/n},
	\end{equation*}
	so that the second formula, for all $k$, amounts to \eqref{eq:yamamoto} for compact operators. This is how \eqref{eq:yamamoto} is derived in \cite{Davis}, and conversely the second formula follows from any proof of \eqref{eq:yamamoto} that does not pass through antisymmetric tensor powers, such as the one in \cite{NS}. The direct proof given above keeps the present paper independent of \eqref{eq:yamamoto}.
	
	The last lemma of this section is needed in order to identify the polar decomposition of $\wedge^kT$.
	
	\begin{lemma}\label{lem:ranges}
		Let $k\in\mathbb{N}$ with $k\le d$.
		\begin{enumerate}[\rm($i$)]
			\item If $Q$ is the orthogonal projection of $\mathcal{H}$ onto a closed subspace $\mathcal{L}$, then $\wedge^kQ$ is the orthogonal projection of $\wedge^k\mathcal{H}$ onto $\wedge^k\mathcal{L}$.
			\item For every $T\in \mathfrak{B}(\mathcal{H})$,  $\overline{\operatorname{ran}(\wedge^kT)}=\wedge^k(\overline{\operatorname{ran}T})$.
		\end{enumerate}
	\end{lemma}
	
	\begin{proof}
		($i$) By Lemma \ref{lem:wedge}($i$) we have $(\wedge^kQ)^2=\wedge^k(Q^2)=\wedge^kQ$ and $(\wedge^kQ)^*=\wedge^k(Q^*)=\wedge^kQ$, so $\wedge^kQ$ is a self-adjoint idempotent and hence the orthogonal projection onto its range. It remains to identify that range.
		
		The vectors $y_1\wedge\cdots\wedge y_k$ with $y_1,\ldots,y_k\in\mathcal{H}$ span a dense subspace of $\wedge^k\mathcal{H}$, and $\wedge^kQ$ maps each of them to $Qy_1\wedge\cdots\wedge Qy_k$, which lies in $\wedge^k\mathcal{L}$. Since $\wedge^kQ$ is bounded and $\wedge^k\mathcal{L}$ is closed, the range of $\wedge^kQ$ is contained in $\wedge^k\mathcal{L}$. Conversely, let $x_1,\ldots,x_k\in\mathcal{L}$. Then $Qx_i=x_i$ for every $i$, and therefore
		\begin{equation*}
			\wedge^kQ\,(x_1\wedge\cdots\wedge x_k)=Qx_1\wedge\cdots\wedge Qx_k=x_1\wedge\cdots\wedge x_k.
		\end{equation*}
		Every vector of this form thus lies in the range of $\wedge^kQ$, and that range is a closed subspace, so it contains their closed linear span $\wedge^k\mathcal{L}$.
		
		($ii$) Put $\mathcal{L}=\overline{\operatorname{ran}T}$. For $y_1,\ldots,y_k\in\mathcal{H}$ we have
		\begin{equation*}
			\wedge^kT\,(y_1\wedge\cdots\wedge y_k)=Ty_1\wedge\cdots\wedge Ty_k\in\wedge^k\mathcal{L}.
		\end{equation*}
		The argument used in ($i$) for $\wedge^kQ$ therefore shows that the range of $\wedge^kT$ is contained in $\wedge^k\mathcal{L}$. Since $\wedge^k\mathcal{L}$ is closed, it also contains the closure of this range.
		
		For the reverse inclusion we first record that the wedge product is jointly continuous. It is multilinear, and by Hadamard's inequality for Gram determinants
		\begin{equation*}
			\|x_1\wedge\cdots\wedge x_k\|^2=\det\big(\langle x_i,x_j\rangle\big)_{i,j=1}^{k}\le\|x_1\|^2\cdots\|x_k\|^2,
		\end{equation*}
		so it is a bounded multilinear map and therefore continuous.
		
		Let now $x_1,\ldots,x_k\in\mathcal{L}$ and choose, for each $i$, a sequence $(y_i^{(m)})_{m\ge1}$ in $\mathcal{H}$ with $Ty_i^{(m)}\to x_i$ as $m\to\infty$. By the continuity just established,
		\begin{equation*}
			\wedge^kT\,\big(y_1^{(m)}\wedge\cdots\wedge y_k^{(m)}\big)=Ty_1^{(m)}\wedge\cdots\wedge Ty_k^{(m)}\longrightarrow x_1\wedge\cdots\wedge x_k,
		\end{equation*}
		and for every $m$ the vector on the left belongs to the range of $\wedge^kT$. Hence $x_1\wedge\cdots\wedge x_k$ belongs to the closure of that range, which is a closed subspace and therefore contains the closed linear span of all such vectors, namely $\wedge^k\mathcal{L}$.
	\end{proof}

	\bigskip 
	\section{The Aluthge transform and \texorpdfstring{$\wedge^k$}{the antisymmetric tensor powers}}\label{sec:commute}
	\bigskip
	
	We now relate $\wedge^k$ to the generalized Aluthge transforms. The first lemma is the identity \eqref{eq:commute} from the introduction, in a slightly more precise form.
	
	\begin{lemma}\label{lem:commute}
		Let $T\in \mathfrak{B}(\mathcal{H})$ have the polar decomposition $T=U|T|$ and let $k\in\mathbb{N}$ with $k\le d$. Then
		\begin{equation*}
			\wedge^kT=(\wedge^kU)(\wedge^k|T|)
		\end{equation*}
		is the polar decomposition of $\wedge^kT$. In particular $|\wedge^kT|=\wedge^k|T|$, and
		\begin{equation*}
			\Delta_t(\wedge^kT)=\wedge^k\Delta_t(T),\qquad 0<t<1.
		\end{equation*}
		Consequently $\Delta_t^n(\wedge^kT)=\wedge^k\Delta_t^n(T)$ for every $n\ge0$.
	\end{lemma}
	
	\begin{proof}
		By Lemma \ref{lem:wedge}($i$) and ($iii$), $\wedge^k|T|$ is positive and $$(\wedge^k|T|)^2=\wedge^k(T^*T)=(\wedge^kT)^*(\wedge^kT),$$ so that $\wedge^k|T|=|\wedge^kT|$. Put $V=\wedge^kU$. Then $V^*V=\wedge^k(U^*U)$, and $U^*U$ is the orthogonal projection onto the initial space $(\ker T)^\perp=\overline{\operatorname{ran}|T|}$ of $U$. By Lemma \ref{lem:ranges}, $V^*V$ is the orthogonal projection onto
		\begin{equation*}
			\wedge^k\big(\overline{\operatorname{ran}|T|}\big)=\overline{\operatorname{ran}(\wedge^k|T|)}=\overline{\operatorname{ran}|\wedge^kT|}=(\ker\wedge^kT)^\perp.
		\end{equation*}
		Thus $V$ is a partial isometry with initial space $(\ker\wedge^kT)^\perp$, and $$\wedge^kT=\wedge^k(U|T|)=V|\wedge^kT|$$ by Lemma \ref{lem:wedge}($i$). By the uniqueness of the polar decomposition, $V$ is the partial isometry in the polar decomposition of $\wedge^kT$. Let now $0<t<1$. By Lemma \ref{lem:power} we have $$\wedge^k(|T|^{t})=(\wedge^k|T|)^{t}=|\wedge^kT|^{t},$$ and the same identity holds with $t$ replaced by $1-t$, so that Lemma \ref{lem:wedge}($i$) gives
		\begin{align*}
			\wedge^k\Delta_t(T)&=\wedge^k\big(|T|^{t}\big)\,\wedge^kU\,\wedge^k\big(|T|^{1-t}\big)\\&=|\wedge^kT|^{t}\,V\,|\wedge^kT|^{1-t}=\Delta_t(\wedge^kT).
		\end{align*}
		The last assertion follows by induction on $n$.
	\end{proof}

	The remaining results of this section concern the Aluthge transform of a compact operator. The following fact is well known and is needed to show that this transform is again compact.
	
	\begin{lemma}\label{lem:compact-power}
		Let $T\in\mathfrak{B}(\mathcal{H})$ be positive and compact and let $\alpha>0$. Then $T^{\alpha}$ is compact.
	\end{lemma}
	
	\begin{proof}
		A bounded operator on a Hilbert space is compact if and only if it maps weakly convergent sequences to norm convergent sequences \cite[Theorem VI.11]{RS}. Let $(x_n)_{n\ge1}$ be a sequence in $\mathcal{H}$ converging weakly to $0$. It is bounded, and $\|Tx_n\|\to0$ because $T$ is compact. Hence
		\begin{equation*}
			\|T^{1/2}x_n\|^2=\langle T^{1/2}x_n,T^{1/2}x_n\rangle=\langle Tx_n,x_n\rangle\le\|Tx_n\|\,\|x_n\|\longrightarrow0,
		\end{equation*}
		so $T^{1/2}$ is compact.
		
		Since $T^{1/2}$ is again positive, the same argument may be repeated, and induction on $m$ gives that $T^{2^{-m}}$ is compact for every integer $m\ge0$.
		
		Let now $\alpha>0$ be arbitrary and choose an integer $m\ge0$ with $2^{-m}<\alpha$. Then
		\begin{equation*}
			T^{\alpha}=T^{\alpha-2^{-m}}T^{2^{-m}},
		\end{equation*}
		where $T^{\alpha-2^{-m}}$ is bounded and $T^{2^{-m}}$ is compact. Hence $T^{\alpha}$ is compact.
	\end{proof}
	
	The next lemma shows that $T$ and $\Delta_t(T)$ have the same non-zero eigenvalues, with the same algebraic multiplicities. In order to state this as an equality of sequences, the ordering convention of Section \ref{sec:prelim} has to be made precise.
	
	\begin{remark}\label{rem:ordering}
		Eigenvalues of equal modulus may be listed in any order, so the sequence $(\lambda_j(T))_j$ is determined by $T$ only up to a permutation within each group of equal moduli. Whenever two compact operators have the same non-zero eigenvalues with the same algebraic multiplicities, as in Lemma \ref{lem:aluthge-basic} below, we list them in the same order for both, so that the two sequences coincide. The results of Sections \ref{sec:main}--\ref{sec:duggal} do not depend on this convention, since only the moduli $|\lambda_j(T)|$ enter them.
	\end{remark}

	\begin{lemma}\label{lem:aluthge-basic}
		Let $T\in\mathcal{K}(\mathcal{H})$ and let $0<t<1$. Then $\Delta_t(T)$ is compact and $\lambda_j(\Delta_t(T))=\lambda_j(T)$ for every $j$.
	\end{lemma}
	
	\begin{proof}
		Put
		\begin{equation*}
			X=U|T|^{1-t},\qquad Y=|T|^{t},
		\end{equation*}
		so that $XY=U|T|^{1-t}|T|^{t}=T$ and $YX=|T|^{t}U|T|^{1-t}=\Delta_t(T)$. The operator $T^*T$ is positive and compact, being the product of a bounded and a compact operator, so Lemma \ref{lem:compact-power}, applied to it with the exponent $t/2$, shows that
		\begin{equation*}
			Y=|T|^{t}=(T^*T)^{t/2}
		\end{equation*}
		is compact. Hence $\Delta_t(T)=YX$ is compact as well.
		
		It remains to compare the eigenvalues. Fix $\mu\neq0$ and $m\in\mathbb{N}$. From $$Y(XY-\mu)=YXY-\mu Y=(YX-\mu)Y$$ we obtain, by induction on $m$,
		\begin{equation*}
			Y(XY-\mu)^m=(YX-\mu)^mY.
		\end{equation*}
		Let $x\in\ker(XY-\mu)^m$. Then $(YX-\mu)^mYx=Y(XY-\mu)^mx=0$, so that $Y$ maps $\ker(XY-\mu)^m$ into $\ker(YX-\mu)^m$. This map is injective. To see it, suppose that $x\in\ker(XY-\mu)^m$ satisfies $Yx=0$. The operators $XY$ and $\mu I$ commute, so the binomial theorem gives
		\begin{equation*}
			(XY-\mu)^mx=\sum_{i=0}^{m}\binom{m}{i}(-\mu)^{m-i}(XY)^{i}x.
		\end{equation*}
		For $i\ge1$ we have $(XY)^{i}x=(XY)^{i-1}X(Yx)=0$, so only the term with $i=0$ survives and $(XY-\mu)^mx=(-\mu)^mx$. The left-hand side vanishes and $\mu\neq0$, whence $x=0$.
		
		Interchanging the roles of $X$ and $Y$ gives $X(YX-\mu)^m=(XY-\mu)^mX$, and the same argument shows that $X$ maps $\ker(YX-\mu)^m$ injectively into $\ker(XY-\mu)^m$. Both kernels are finite-dimensional, because $XY$ and $YX$ are compact and $\mu\neq0$, and the two injections give
		\begin{equation*}
			\dim\ker(XY-\mu)^m=\dim\ker(YX-\mu)^m,\qquad \mu\neq0,\quad m\in\mathbb{N}.
		\end{equation*}
		
		Taking $m=1$ shows that $T=XY$ and $\Delta_t(T)=YX$ have the same non-zero eigenvalues, and taking $m$ large shows that each of them has the same algebraic multiplicity for the two operators. The sequences $(\lambda_j(T))_j$ and $(\lambda_j(\Delta_t(T)))_j$ are therefore obtained by arranging the same multiset of non-zero numbers in non-increasing order of modulus, in the same order by Remark \ref{rem:ordering}, and by continuing with zeros over the same index set, so they coincide.
	\end{proof} 
	
	\bigskip
	\section{Proof of the main theorem}\label{sec:main}
	\bigskip
	
	We begin with the products of the $k$ largest singular values. This is the point where \eqref{eq:commute} and \eqref{eq:yamazaki} enter, and the results of this section and the next are derived from the following theorem.
	
	\begin{theorem}\label{thm:products}
		Let $T\in\mathcal{K}(\mathcal{H})$, let $0<t<1$ and let $k\in\mathbb{N}$ with $k\le d$. Then the sequence
		\begin{equation*}
			P_k(n)=\prod_{j=1}^{k}s_j(\Delta_t^n(T)),\qquad n=0,1,2,\ldots,
		\end{equation*}
		is non-increasing and
		\begin{equation*}
			\lim_{n\to\infty}\prod_{j=1}^{k}s_j(\Delta_t^n(T))=\prod_{j=1}^{k}|\lambda_j(T)|.
		\end{equation*}
	\end{theorem}
	
	\begin{proof}
		By Lemma \ref{lem:aluthge-basic} the operator $\Delta_t^n(T)$ is compact for every $n$, so the first formula of Lemma \ref{lem:wedge-compact}, applied to it, gives
		\begin{equation*}
			P_k(n)=\big\|\wedge^k\Delta_t^n(T)\big\|.
		\end{equation*}
		By Lemma \ref{lem:commute} we have $\wedge^k\Delta_t^n(T)=\Delta_t^n(\wedge^kT)$, and therefore
		\begin{equation}\label{eq:Pk}
			P_k(n)=\big\|\Delta_t^n(\wedge^kT)\big\|,\qquad n=0,1,2,\ldots
		\end{equation} 
		
		We first prove that $P_k$ is non-increasing in $n$. Fix $n$ and put $S=\Delta_t^n(\wedge^kT)$. Applying \eqref{eq:contraction} on $\wedge^k\mathcal{H}$ to $S$ gives $\|\Delta_t(S)\|\le\|S\|$, and $\Delta_t(S)=\Delta_t^{n+1}(\wedge^kT)$, so that \eqref{eq:Pk} yields
		\begin{equation*}
			P_k(n+1)=\big\|\Delta_t^{n+1}(\wedge^kT)\big\|\le\big\|\Delta_t^{n}(\wedge^kT)\big\|=P_k(n).
		\end{equation*}
		
		It remains to identify the limit. By \eqref{eq:Pk}, the sequence $P_k$ is the sequence of norms of the iterated Aluthge transforms of $\wedge^kT$, so \eqref{eq:yamazaki}, applied to $\wedge^kT$ on $\wedge^k\mathcal{H}$, gives
		\begin{equation*}
			\lim_{n\to\infty}P_k(n)=r(\wedge^kT).
		\end{equation*}
		The operator $T$ is compact, so the second formula of Lemma \ref{lem:wedge-compact}, applied to $T$, evaluates this spectral radius and completes the proof:
		\begin{equation*}
			\lim_{n\to\infty}P_k(n)=r(\wedge^kT)=\prod_{j=1}^{k}|\lambda_j(T)|.\qedhere
		\end{equation*}
	\end{proof}
	
	\begin{remark}\label{rem:known}
		The monotonicity of $P_k$ in $n$ is not new. It follows by iterating \eqref{eq:logmaj}, which is stated in \cite[Remark 3.12]{AMS} for every parameter in $[0,1]$, so that it holds for the Duggal transform $\Delta_1$ as well. In finite dimensions the limit is not new either: for $t=1/2$ it is due to Ando \cite{Ando}, and for general $t$ it follows in the same way from the matrix case of \eqref{eq:yamazaki}, see \cite[Theorem 1.4]{HT}. The content of Theorem \ref{thm:products} is the limit for compact operators on infinite-dimensional spaces. The restriction $t<1$ cannot be dropped in any dimension, see Proposition \ref{prop:duggal} and the remark following it.
	\end{remark}
	
	We are now in a position to prove the main result. The individual singular values are obtained from Theorem \ref{thm:products} by taking ratios of consecutive products, and the only point requiring care is the possible vanishing of the denominators.
	
	\begin{proof}[Proof of Theorem \ref{thm:main-intro}]
		The assertion about the products is Theorem \ref{thm:products}, so it remains to prove that $s_j(\Delta_t^n(T))\to|\lambda_j(T)|$ for every $j$. Let $P_k(n)$ be as in Theorem \ref{thm:products}, put $P_0(n)=1$, and let $\rho=\rho(T)$. With this convention the conclusions of Theorem \ref{thm:products} hold for $k=0$ as well, the empty product being $1$. Throughout, $j$ denotes a positive integer.
		
		We first treat the indices $j$ with $j\le\rho$, for which $\lambda_j(T)\neq0$. Fix such a $j$. Then $\lambda_1(T),\ldots,\lambda_{j-1}(T)$ are non-zero as well, so by Theorem \ref{thm:products} the sequence $P_{j-1}(n)$ decreases to $\prod_{i<j}|\lambda_i(T)|>0$. Being non-increasing, it stays above its limit, and in particular $P_{j-1}(n)>0$ for every $n$. Theorem \ref{thm:products} applied to the indices $j$ and $j-1$ therefore gives
		\begin{equation*}
			s_j(\Delta_t^n(T))=\frac{P_j(n)}{P_{j-1}(n)}\longrightarrow\frac{\prod_{i\le j}|\lambda_i(T)|}{\prod_{i<j}|\lambda_i(T)|}=|\lambda_j(T)|.
		\end{equation*}
		
		Suppose that $\rho=\infty$. Every positive integer $j$ then satisfies $j\le\rho$, so the assertion has already been proved for every $j$.
		
		Suppose now that $\rho<\infty$. The indices $j\le\rho$ have been treated, and indices $j>\rho$ occur only when $\rho<d$, so we may assume that $\rho<d$. For such $j$ we have $\lambda_j(T)=0$, and the assertion to be proved is that $s_j(\Delta_t^n(T))\to0$. Since $\rho+1\le d$, Theorem \ref{thm:products} applies to the index $\rho+1$ and gives
		\begin{equation*}
			\lim_{n\to\infty}P_{\rho+1}(n)=\prod_{i\le\rho+1}|\lambda_i(T)|=0,
		\end{equation*}
		because the factor with $i=\rho+1$ vanishes. On the other hand $P_\rho(n)$ decreases to $\prod_{i\le\rho}|\lambda_i(T)|>0$ and therefore stays above this positive number. Consequently
		\begin{equation*}
			s_{\rho+1}(\Delta_t^n(T))=\frac{P_{\rho+1}(n)}{P_\rho(n)}\longrightarrow0,
		\end{equation*}
		and since the singular values are arranged in non-increasing order, $$s_j(\Delta_t^n(T))\le s_{\rho+1}(\Delta_t^n(T))$$ for every $j\ge\rho+1$. Hence $s_j(\Delta_t^n(T))\to0=|\lambda_j(T)|$ for all these indices as well.
	\end{proof}
	
	\bigskip
	\section{Unitarily invariant norms and Schatten classes}\label{sec:norms}
	\bigskip
	
	We now pass from products to sums, which requires a word on majorization. For $x,y\in\mathbb{R}^m$ with non-negative coordinates arranged in non-increasing order we write $x\prec_{w\log}y$ if $\prod_{j\le k}x_j\le\prod_{j\le k}y_j$ for $k=1,\ldots,m$, and $x\prec_w y$ if $\sum_{j\le k}x_j\le\sum_{j\le k}y_j$ for $k=1,\ldots,m$. The same notation is used for infinite sequences, with the inequalities required for every $k$. Since $u\mapsto e^{u}$ is convex and increasing, $x\prec_{w\log}y$ implies $x\prec_w y$, see \cite[Corollary 1.10]{Simon} or \cite[Chapter II.3]{Bhatia}. For a compact operator $T$ we write $s(T)=(s_j(T))_j$ and $\lambda(T)=(\lambda_j(T))_j$. Weyl's inequality $|\lambda(T)|\prec_{w\log}s(T)$, see \cite[Theorem 1.14]{Simon}, thus yields Weyl's majorant theorem $|\lambda(T)|\prec_{w}s(T)$, see \cite[Theorem 1.15]{Simon}. For $k\in\mathbb{N}$ with $k\le d$ the Ky Fan norm of $T\in\mathcal{K}(\mathcal{H})$ is
	\begin{equation*}
		\|T\|_{(k)}=\sum_{j=1}^{k}s_j(T).
	\end{equation*}
	For $k=1$ this is the operator norm.
	
	\begin{corollary}\label{cor:kyfan}
		Let $T\in\mathcal{K}(\mathcal{H})$, let $0<t<1$ and let $k\in\mathbb{N}$ with $k\le d$. Then the sequence $\|\Delta_t^n(T)\|_{(k)}$, $n=0,1,2,\ldots$, is non-increasing, and
		\begin{equation*}
			\lim_{n\to\infty}\|\Delta_t^n(T)\|_{(k)}=\sum_{j=1}^{k}|\lambda_j(T)|.
		\end{equation*}
	\end{corollary}
	
	\begin{proof}
		Theorem \ref{thm:products}, applied for all indices up to $k$, says that the vector of the first $k$ singular values of $\Delta_t^{n+1}(T)$ is weakly log-majorized by that of $\Delta_t^n(T)$. Hence it is weakly majorized by it, which gives $\|\Delta_t^{n+1}(T)\|_{(k)}\le\|\Delta_t^n(T)\|_{(k)}$. The limit follows from Theorem \ref{thm:main-intro}, since the sum is finite.
	\end{proof}
	
	In finite dimensions, where Theorem \ref{thm:main-intro} is known, the same argument covers all unitarily invariant norms. We include the short derivation because it shows exactly what is missing in infinite dimensions. A norm $\Phi$ on $\mathbb{R}^d$ is called a symmetric gauge function if
	\begin{equation*}
		\Phi(\varepsilon_1x_{\pi(1)},\ldots,\varepsilon_dx_{\pi(d)})=\Phi(x_1,\ldots,x_d)
	\end{equation*}
	for every permutation $\pi$ of $\{1,\ldots,d\}$ and all signs $\varepsilon_1,\ldots,\varepsilon_d\in\{-1,1\}$. The basic examples are the $\ell^p$ norms, $1\le p\le\infty$. A norm $|\mkern-2mu|\mkern-2mu|\cdot|\mkern-2mu|\mkern-2mu|$ on the algebra of $d\times d$ complex matrices is unitarily invariant if $|\mkern-2mu|\mkern-2mu|V_1TV_2|\mkern-2mu|\mkern-2mu|=|\mkern-2mu|\mkern-2mu|T|\mkern-2mu|\mkern-2mu|$ for all unitary matrices $V_1$ and $V_2$. By a classical theorem of von Neumann, the two notions are in a bijective correspondence: every symmetric gauge function $\Phi$ defines a unitarily invariant norm through
	\begin{equation*}
		|\mkern-2mu|\mkern-2mu|T|\mkern-2mu|\mkern-2mu|=\Phi\big(s_1(T),\ldots,s_d(T)\big),
	\end{equation*}
	and every unitarily invariant norm arises in this way from exactly one $\Phi$. Under this correspondence the $\ell^p$ norms produce the Schatten $p$-norms $\|T\|_p=\big(\sum_{j=1}^{d}s_j(T)^{p}\big)^{1/p}$, and $\ell^\infty$ produces the operator norm. The property we need is Fan's dominance principle: if $x$ and $y$ have non-negative coordinates arranged in non-increasing order and $x\prec_w y$, then $\Phi(x)\le\Phi(y)$. For all these facts see \cite[Chapter IV]{Bhatia}.
	
	\begin{corollary}\label{cor:uinorm}
		Let $T$ be a $d\times d$ complex matrix, let $0<t<1$ and let $|\mkern-2mu|\mkern-2mu|\cdot|\mkern-2mu|\mkern-2mu|$ be a unitarily invariant norm with symmetric gauge function $\Phi$. Then the sequence $|\mkern-2mu|\mkern-2mu|\Delta_t^n(T)|\mkern-2mu|\mkern-2mu|$ is non-increasing and
		\begin{equation*}
			\lim_{n\to\infty}|\mkern-2mu|\mkern-2mu|\Delta_t^n(T)|\mkern-2mu|\mkern-2mu|=\Phi\big(|\lambda_1(T)|,\ldots,|\lambda_d(T)|\big).
		\end{equation*}
		In particular, $\lim_{n\to\infty}\|\Delta_t^n(T)\|_p=\big(\sum_{j=1}^{d}|\lambda_j(T)|^p\big)^{1/p}$ for every Schatten $p$-norm, $1\le p<\infty$.
	\end{corollary}
	
	\begin{proof}
		Theorem \ref{thm:products} says that
		\begin{equation*}
			s(\Delta_t^{n+1}(T))\prec_{w\log}s(\Delta_t^{n}(T)),\qquad\text{hence}\qquad s(\Delta_t^{n+1}(T))\prec_{w}s(\Delta_t^{n}(T)),
		\end{equation*}
		weak log-majorization implying weak majorization. Fan's dominance principle then gives
		\begin{equation*}
			|\mkern-2mu|\mkern-2mu|\Delta_t^{n+1}(T)|\mkern-2mu|\mkern-2mu|\le|\mkern-2mu|\mkern-2mu|\Delta_t^{n}(T)|\mkern-2mu|\mkern-2mu|.
		\end{equation*} By Theorem \ref{thm:main-intro} the vector $s(\Delta_t^{n}(T))$ converges to $\big(|\lambda_1(T)|,\ldots,|\lambda_d(T)|\big)$ in each coordinate, hence in $\mathbb{R}^d$, and $\Phi$ is continuous, being a norm on $\mathbb{R}^d$. Therefore
		\begin{equation*}
			|\mkern-2mu|\mkern-2mu|\Delta_t^{n}(T)|\mkern-2mu|\mkern-2mu|=\Phi\big(s(\Delta_t^{n}(T))\big)\longrightarrow\Phi\big(|\lambda_1(T)|,\ldots,|\lambda_d(T)|\big).\qedhere
		\end{equation*}
	\end{proof}
	
	\begin{remark}
		Corollary \ref{cor:uinorm} is also a consequence of the theorem of Antezana, Pujals and Stojanoff, according to which the iterates $\Delta_t^n(T)$ of every matrix $T$ converge to a normal matrix with the same eigenvalues for every $t\in(0,1)$, see \cite[Theorem 6.1]{APS}. The diagonalizable case was settled earlier in \cite{APS2}. The proof given above does not use that theorem. In infinite dimensions no such convergence theorem is available. The iterates of a bounded operator need not converge even in the weak operator topology \cite{CJL}, and for compact operators the convergence of the iterates is open as far as we know, see \cite{JKP2} and \cite[Section 3.3]{APS3}. For results on special classes of operators see \cite{OY}. Theorem \ref{thm:main-intro} and Corollary \ref{cor:kyfan} show that their singular values and Ky Fan norms converge regardless.
	\end{remark}
	
	We now turn to compact operators on an infinite-dimensional space, and to unitarily invariant norms defined on the corresponding norm ideals of $\mathcal{K}(\mathcal{H})$, see \cite[Chapter III]{GK}. The monotonicity extends at once. Fan's dominance property is available in this setting in a strong form: if $Y$ belongs to the ideal $\mathcal{J}_{|\mkern-2mu|\mkern-2mu|\cdot|\mkern-2mu|\mkern-2mu|}$ and $X$ is compact with $\|X\|_{(k)}\le\|Y\|_{(k)}$ for every $k$, then $X$ belongs to that ideal as well and $|\mkern-2mu|\mkern-2mu|X|\mkern-2mu|\mkern-2mu|\le|\mkern-2mu|\mkern-2mu|Y|\mkern-2mu|\mkern-2mu|$, see \cite[Chapter III]{GK} and \cite{JK}. Applying this with $X=\Delta_t^{n+1}(T)$ and $Y=\Delta_t^{n}(T)$, Corollary \ref{cor:kyfan} gives the following.
	
	\begin{corollary}\label{cor:ideal-monotone}
		Let $T\in\mathcal{K}(\mathcal{H})$, let $0<t<1$ and let $|\mkern-2mu|\mkern-2mu|\cdot|\mkern-2mu|\mkern-2mu|$ be a unitarily invariant norm with norm ideal $\mathcal{J}_{|\mkern-2mu|\mkern-2mu|\cdot|\mkern-2mu|\mkern-2mu|}$. If $T\in\mathcal{J}_{|\mkern-2mu|\mkern-2mu|\cdot|\mkern-2mu|\mkern-2mu|}$, then $\Delta_t^{n}(T)\in\mathcal{J}_{|\mkern-2mu|\mkern-2mu|\cdot|\mkern-2mu|\mkern-2mu|}$ for every $n$, and the sequence $|\mkern-2mu|\mkern-2mu|\Delta_t^{n}(T)|\mkern-2mu|\mkern-2mu|$ is non-increasing, hence convergent.
	\end{corollary}
	
	The limit is a more delicate matter, and here the passage from Ky Fan norms to a general unitarily invariant norm is not automatic. Fan's dominance property compares two fixed operators, whereas the limit requires passing from the convergence of $\|\Delta_t^{n}(T)\|_{(k)}$ for each fixed $k$ to the convergence of $|\mkern-2mu|\mkern-2mu|\Delta_t^{n}(T)|\mkern-2mu|\mkern-2mu|$, and no such implication holds. Let $d=\infty$ and let $X_n$ be the operator which is diagonal with respect to a fixed orthonormal basis of $\mathcal{H}$, with the first $n$ diagonal entries equal to $1/n$ and all others equal to $0$. Then $s_j(X_n)=1/n$ for $j\le n$ and $s_j(X_n)=0$ for $j>n$, so that $\|X_n\|_{(k)}=k/n$ for $n\ge k$, which tends to $0$ for every fixed $k$, while the trace norm $\|X_n\|_1=\sum_{j}s_j(X_n)=\sup_k\|X_n\|_{(k)}$ equals $1$ for every $n$. In finite dimensions this cannot happen, since $\|\cdot\|_{(d)}$ is the largest Ky Fan norm and the singular value vector converges in $\mathbb{R}^d$, which is how the limit in Corollary \ref{cor:uinorm} was obtained. For $d=\infty$ the supremum over $k$ need not commute with the limit in $n$, as the example shows. Thus the convergence of all Ky Fan norms in Corollary \ref{cor:kyfan} does not by itself identify the limit of $|\mkern-2mu|\mkern-2mu|\Delta_t^{n}(T)|\mkern-2mu|\mkern-2mu|$, which exists by Corollary \ref{cor:ideal-monotone}.
	
	What restores the passage to the limit is a bound on the tails $\sum_{j>k}s_j(\Delta_t^{n}(T))^{p}$ that is uniform in $n$, and the log-majorization provides one as soon as $T$ lies in a Schatten class strictly smaller than the one in which the limit is taken. For $0<p<\infty$ we write $\|X\|_p=\big(\sum_{j}s_j(X)^{p}\big)^{1/p}$ and $\mathcal{S}_p=\{X\in\mathcal{K}(\mathcal{H}):\|X\|_p<\infty\}$. For $p\ge1$ this is the Schatten $p$-norm and $\mathcal{S}_p$ is its norm ideal, whereas for $0<p<1$ the functional $\|\cdot\|_p$ is only a quasi-norm; the statements below hold in both cases. For $p=1$ the ideal $\mathcal{S}_1$ is the trace class, $\|\cdot\|_1$ is the trace norm, and $\operatorname{tr}$ denotes the trace on $\mathcal{S}_1$.
	
	\begin{theorem}\label{thm:schatten}
		Let $0<q<p<\infty$, let $T\in\mathcal{S}_q$ and let $0<t<1$. Then $\Delta_t^{n}(T)\in\mathcal{S}_q$ for every $n$, the sequence $\|\Delta_t^{n}(T)\|_p$ is non-increasing, and
		\begin{equation*}
			\lim_{n\to\infty}\|\Delta_t^{n}(T)\|_p=\Big(\sum_{j\ge1}|\lambda_j(T)|^p\Big)^{1/p}.
		\end{equation*}
	\end{theorem}
	
	\begin{proof}
		Write $s_j^{(n)}=s_j(\Delta_t^{n}(T))$ and $\lambda_j=\lambda_j(T)$. By Theorem \ref{thm:products} the sequence $(s_j^{(n)})_j$ is weakly log-majorized by $(s_j(T))_j$ for every $n$, and $(s_j^{(n+1)})_j$ is weakly log-majorized by $(s_j^{(n)})_j$. For every $\alpha>0$ the function $u\mapsto u^{\alpha}$ is increasing on $[0,\infty)$ and $u\mapsto e^{\alpha u}$ is convex, so weak log-majorization of two finite sequences is turned by this function into weak majorization \cite[Corollary 1.10]{Simon}. Applying this with $\alpha=q$ to the first of the two log-majorizations and with $\alpha=p$ to the second, in each case to the first $m$ terms of the sequences, and letting $m\to\infty$, we obtain
		\begin{equation}\label{eq:tail}
			\sum_{j\ge1}\big(s_j^{(n)}\big)^{q}\le\sum_{j\ge1}s_j(T)^{q}=\|T\|_q^{q},\qquad n=0,1,2,\ldots,
		\end{equation}
		and
		\begin{equation*}
			\sum_{j\ge1}\big(s_j^{(n+1)}\big)^{p}\le\sum_{j\ge1}\big(s_j^{(n)}\big)^{p},\qquad n=0,1,2,\ldots
		\end{equation*}
		The first inequality shows that every $\Delta_t^{n}(T)$ lies in $\mathcal{S}_q$, and the second one is the monotonicity of $\|\Delta_t^{n}(T)\|_p$.
		
		Since the singular values are arranged in non-increasing order, \eqref{eq:tail} gives
		\begin{equation*}
			(k+1)\big(s_{k+1}^{(n)}\big)^{q}\le\|T\|_q^{q},\qquad\text{that is,}\qquad s_{k+1}^{(n)}\le\frac{\|T\|_q}{(k+1)^{1/q}},\qquad k\ge0,\quad n\ge0,
		\end{equation*}
		a bound which does not depend on $n$. Using it together with \eqref{eq:tail} and $p>q$, we estimate the tails uniformly in $n$:
		\begin{equation*}
			\sum_{j>k}\big(s_j^{(n)}\big)^{p}\le\big(s_{k+1}^{(n)}\big)^{p-q}\sum_{j>k}\big(s_j^{(n)}\big)^{q}\le\frac{\|T\|_q^{\,p-q}}{(k+1)^{(p-q)/q}}\,\|T\|_q^{q}=\frac{\|T\|_q^{\,p}}{(k+1)^{(p-q)/q}}.
		\end{equation*}
		Denote the right-hand side by $\eta_k$, so that $\eta_k\to0$ as $k\to\infty$. By Theorem \ref{thm:main-intro} we have $s_j^{(n)}\to|\lambda_j|$ as $n\to\infty$ for every $j$, so Fatou's lemma, applied to the counting measure on $\{j:j>k\}$, gives
		\begin{equation*}
			\sum_{j>k}|\lambda_j|^{p}\le\liminf_{n\to\infty}\sum_{j>k}\big(s_j^{(n)}\big)^{p}\le\eta_k,
		\end{equation*}
		and in particular $\sum_{j\ge1}|\lambda_j|^{p}$ is finite.
		
		Let $k$ be fixed. Splitting both sums at $k$ we have
		\begin{equation*}
			\|\Delta_t^{n}(T)\|_p^{p}-\sum_{j\ge1}|\lambda_j|^{p}=\Big(\sum_{j\le k}\big(s_j^{(n)}\big)^{p}-\sum_{j\le k}|\lambda_j|^{p}\Big)+\sum_{j>k}\big(s_j^{(n)}\big)^{p}-\sum_{j>k}|\lambda_j|^{p},
		\end{equation*}
		and both tails are at most $\eta_k$, so that
		\begin{equation*}
			\Big|\,\|\Delta_t^{n}(T)\|_p^{p}-\sum_{j\ge1}|\lambda_j|^{p}\Big|\le\Big|\sum_{j\le k}\big(s_j^{(n)}\big)^{p}-\sum_{j\le k}|\lambda_j|^{p}\Big|+2\eta_k,\qquad n\ge0.
		\end{equation*}
		The sum over $j\le k$ is finite, and $s_j^{(n)}\to|\lambda_j|$ for each $j$ by Theorem \ref{thm:main-intro}, so the first term on the right tends to $0$ as $n\to\infty$. Hence
		\begin{equation*}
			\limsup_{n\to\infty}\Big|\,\|\Delta_t^{n}(T)\|_p^{p}-\sum_{j\ge1}|\lambda_j|^{p}\Big|\le2\eta_k
		\end{equation*}
		for every $k$. Since $\eta_k\to0$ as $k\to\infty$, the limit superior is $0$, that is, $\|\Delta_t^{n}(T)\|_p^{p}\to\sum_{j\ge1}|\lambda_j|^{p}$. Taking $p$-th roots gives the assertion.
	\end{proof}
	
	For trace class operators Theorem \ref{thm:schatten} gives the following.
	
	\begin{corollary}\label{cor:trace-class}
		Let $T$ be of trace class and let $0<t<1$. Then
		\begin{equation*}
			\lim_{n\to\infty}\|\Delta_t^{n}(T)\|_p=\Big(\sum_{j\ge1}|\lambda_j(T)|^{p}\Big)^{1/p}
		\end{equation*}
		for every $p>1$.
	\end{corollary}
	
	\begin{proof}
		Apply Theorem \ref{thm:schatten} with $q=1$.
	\end{proof}
	
	\begin{remark}\label{rem:ideals}
		Theorem \ref{thm:schatten} leaves out exactly the diagonal case $p=q$, and the smallest instance of it is the trace norm of a trace class operator that lies in no $\mathcal{S}_q$ with $q<1$. There the tail estimate above is unavailable, and we do not know whether
		\begin{equation*}
			\lim_{n\to\infty}\|\Delta_t^{n}(T)\|_1=\sum_{j\ge1}|\lambda_j(T)|
		\end{equation*}
		always holds. Two things are known in this case. The limit on the left exists, by Corollary \ref{cor:ideal-monotone}, and it dominates the right-hand side, since for every $k$
		\begin{equation*}
			\sum_{j\le k}|\lambda_j(T)|=\lim_{n\to\infty}\|\Delta_t^{n}(T)\|_{(k)}\le\lim_{n\to\infty}\|\Delta_t^{n}(T)\|_{1}
		\end{equation*}
		by Corollary \ref{cor:kyfan}, and the right-hand side is finite by Weyl's majorant theorem. Moreover, the trace itself is preserved. Indeed, $\Delta_t^{n}(T)$ is of trace class by Corollary \ref{cor:ideal-monotone} and has the same non-zero eigenvalues as $T$, with the same algebraic multiplicities, by repeated application of Lemma \ref{lem:aluthge-basic}, so Lidskii's theorem \cite[Chapter 3]{Simon} gives
		\begin{equation*}
			\operatorname{tr}\Delta_t^{n}(T)=\sum_{j\ge1}\lambda_j(T)=\operatorname{tr}T,\qquad n=0,1,2,\ldots
		\end{equation*}
	\end{remark}
	
	In the diagonal case the limit can nevertheless be identified for some classes of operators. We begin with weighted shifts, whose iterates tend to zero in every Schatten class containing them. Recall that an operator $S$ on an infinite-dimensional Hilbert space is a unilateral weighted shift if there are an orthonormal basis $(e_k)_{k\ge1}$ of the space and complex numbers $w_k$, the weights of $S$, such that $Se_k=w_ke_{k+1}$ for every $k$. We refer to \cite{Shields} for the general theory of these operators. Since $S^*Se_k=|w_k|^{2}e_k$, the operator $|S|$ is diagonal with the entries $|w_k|$, so that $S\in\mathcal{S}_p$ if and only if $\sum_k|w_k|^{p}<\infty$, in which case $\|S\|_p^{p}=\sum_k|w_k|^{p}$. A unilateral weighted shift has no non-zero eigenvalues, see \cite[Theorem 8]{Shields} for the case of non-zero weights. In general, if $Sx=\mu x$ with $\mu\neq0$ and $x=\sum_kx_ke_k$, then comparing the coefficients of $e_1$ and of $e_{k+1}$ on both sides gives $\mu x_1=0$ and $\mu x_{k+1}=w_kx_k$ for every $k$, so that $x=0$.
	
	\begin{proposition}\label{prop:shift}
		Let $0<p<\infty$, let $0<t<1$, and let $S\in\mathcal{S}_p$ be a unilateral weighted shift. Then
		\begin{equation*}
			\lim_{n\to\infty}\|\Delta_t^{n}(S)\|_p=0.
		\end{equation*}
	\end{proposition}
	
	\begin{proof}
		Let $(e_k)_{k\ge1}$ be an orthonormal basis and $(w_k)_{k\ge1}$ weights with $Se_k=w_ke_{k+1}$ for every $k$. We first reduce to non-negative weights. Write $w_k=|w_k|\theta_k$ with $|\theta_k|=1$, taking $\theta_k=1$ when $w_k=0$, and let $D$ be the diagonal operator with $De_k=\alpha_ke_k$, where $\alpha_1=1$ and $\alpha_{k+1}=\theta_k\alpha_k$. Since $|\alpha_k|=1$ for every $k$, the operator $D$ is unitary, and
		\begin{equation*}
			D^*SDe_k=\alpha_kD^*Se_k=\alpha_kw_kD^*e_{k+1}=\alpha_kw_k\overline{\alpha_{k+1}}\,e_{k+1}=\alpha_k|w_k|\theta_k\overline{\theta_k}\,\overline{\alpha_k}\,e_{k+1}=|w_k|e_{k+1}.
		\end{equation*}
		Thus $S'=D^*SD$ is the weighted shift with weights $|w_k|$, and $S=DS'D^*$. If $S'=U'|S'|$ is the polar decomposition of $S'$, then $D|S'|D^*$ is positive with $(D|S'|D^*)^2=DS'^*S'D^*=S^*S$, so that $D|S'|D^*=|S|$, and $DU'D^*$ is a partial isometry with initial space $D(\ker S')^\perp=(\ker S)^\perp$. By the uniqueness of the polar decomposition, $S=(DU'D^*)|S|$ is the polar decomposition of $S$, and the functional calculus gives $|S|^{\alpha}=D|S'|^{\alpha}D^*$ for every $\alpha>0$. Hence $\Delta_t(S)=D\Delta_t(S')D^*$, and by induction $\Delta_t^{n}(S)=D\Delta_t^{n}(S')D^*$ for every $n$. Unitarily equivalent operators have the same singular values, so $\|\Delta_t^{n}(S)\|_p=\|\Delta_t^{n}(S')\|_p$ for every $n$. We may therefore assume that $w_k\ge0$ for every $k$, so that $\|S\|_p^{p}=\sum_kw_k^{p}$.
		
		Now $|S|e_k=w_ke_k$, and the partial isometry $U$ in the polar decomposition of $S$ maps $e_k$ to $e_{k+1}$ when $w_k>0$ and to $0$ otherwise. Consequently
		\begin{equation*}
			\Delta_t(S)e_k=|S|^{t}U|S|^{1-t}e_k=w_k^{1-t}|S|^{t}Ue_k=w_k^{1-t}w_{k+1}^{t}e_{k+1},
		\end{equation*}
		that is, $\Delta_t(S)$ is the weighted shift with the non-negative weights $w_k^{1-t}w_{k+1}^{t}$. We claim that $\Delta_t^{n}(S)$ is the weighted shift with weights
		\begin{equation}\label{eq:weights}
			w_k^{(n)}=\prod_{i=0}^{n}w_{k+i}^{\,b_{n,i}},\qquad b_{n,i}=\binom{n}{i}t^{i}(1-t)^{n-i}.
		\end{equation}
		For $n=0$ this is the definition of $S$, since $b_{0,0}=1$. If it holds for some $n$, then applying the formula for $\Delta_t(S)$ to the weighted shift $\Delta_t^{n}(S)$ shows that $\Delta_t^{n+1}(S)$ is the weighted shift with weights
		\begin{equation*}
			w_k^{(n+1)}=\big(w_k^{(n)}\big)^{1-t}\big(w_{k+1}^{(n)}\big)^{t}=\prod_{i=0}^{n}w_{k+i}^{(1-t)b_{n,i}}\prod_{i=0}^{n}w_{k+1+i}^{\,t\,b_{n,i}}=\prod_{m=0}^{n+1}w_{k+m}^{(1-t)b_{n,m}+t\,b_{n,m-1}},
		\end{equation*}
		where $b_{n,-1}=b_{n,n+1}=0$, and the exponents are
		\begin{align*}
			(1-t)b_{n,m}+t\,b_{n,m-1}&=\Big[\binom{n}{m}+\binom{n}{m-1}\Big]t^{m}(1-t)^{n+1-m}\\&=\binom{n+1}{m}t^{m}(1-t)^{n+1-m}=b_{n+1,m}.
		\end{align*}
		This proves \eqref{eq:weights}. All exponents $b_{n,i}$ are positive, so a zero weight among $w_k,\ldots,w_{k+n}$ gives $w_k^{(n)}=0$, in accordance with the formula for one step.
		
		By \eqref{eq:weights}, $|\Delta_t^{n}(S)|$ is the diagonal operator with the entries $w_k^{(n)}$. Since $\sum_{i}b_{n,i}=1$, the inequality between weighted geometric and arithmetic means gives
		\begin{equation*}
			\big(w_k^{(n)}\big)^{p}=\prod_{i=0}^{n}\big(w_{k+i}^{p}\big)^{b_{n,i}}\le\sum_{i=0}^{n}b_{n,i}w_{k+i}^{p},
		\end{equation*}
		and therefore
		\begin{equation*}
			\|\Delta_t^{n}(S)\|_p^{p}=\sum_{k\ge1}(w_k^{(n)})^{p}\le\sum_{i=0}^{n}b_{n,i}\sum_{k\ge1}w_{k+i}^{p}=\sum_{i=0}^{n}b_{n,i}R_i,\qquad R_i=\sum_{m>i}w_m^{p}.
		\end{equation*}
		The tails $R_i$ of the convergent series $\sum_mw_m^{p}=\|S\|_p^{p}$ decrease to $0$. For fixed $i$ we have $b_{n,i}\le n^{i}(1-t)^{n-i}\to0$ as $n\to\infty$, because $0<1-t<1$, and hence $\sum_{i<M}b_{n,i}\to0$ for every fixed $M$. Splitting the sum at $M$ and using $R_i\le R_0$ for $i<M$, $R_i\le R_M$ for $i\ge M$ and $\sum_ib_{n,i}=1$, we obtain
		\begin{equation*}
			\sum_{i=0}^{n}b_{n,i}R_i\le R_0\sum_{i<M}b_{n,i}+R_M,\qquad n\ge M,
		\end{equation*}
		and the right-hand side tends to $R_M$ as $n\to\infty$. Hence
		\begin{equation*}
			\limsup_{n\to\infty}\|\Delta_t^{n}(S)\|_p^{p}\le R_M
		\end{equation*}
		for every $M$, and letting $M\to\infty$ shows that $\|\Delta_t^{n}(S)\|_p\to0$.
	\end{proof}
	
	Proposition \ref{prop:shift} is not contained in Theorem \ref{thm:schatten}. For instance, the weighted shift with weights $w_k=1/\big((k+1)\log^{2}(k+1)\big)$ is of trace class but belongs to no $\mathcal{S}_q$ with $q<1$, and Proposition \ref{prop:shift} gives $\|\Delta_t^{n}(S)\|_1\to0$ for it. Convergence in the operator norm alone is not new here, since $r(S)=0$ and \eqref{eq:yamazaki} already give $\|\Delta_t^{n}(S)\|\to0$.
	
	\begin{corollary}\label{cor:diagonal}
		Let $0<p<\infty$, let $0<t<1$, and let $T\in\mathcal{S}_p$ be an orthogonal direct sum $T=\bigoplus_{\ell}T_\ell$ in which every $T_\ell$ either acts on a finite-dimensional space or is a unilateral weighted shift. Then
		\begin{equation*}
			\lim_{n\to\infty}\|\Delta_t^{n}(T)\|_p=\Big(\sum_{j\ge1}|\lambda_j(T)|^{p}\Big)^{1/p}.
		\end{equation*}
	\end{corollary}
	
	\begin{proof}
		The polar decomposition and the functional calculus are compatible with orthogonal direct sums, so $\Delta_t(T)=\bigoplus_\ell\Delta_t(T_\ell)$, and by induction $\Delta_t^{n}(T)=\bigoplus_\ell\Delta_t^{n}(T_\ell)$ for every $n$. The singular values of an orthogonal direct sum are those of the summands, and so are the non-zero eigenvalues, counted with algebraic multiplicity. Hence
		\begin{equation*}
			\|\Delta_t^{n}(T)\|_p^{p}=\sum_\ell\|\Delta_t^{n}(T_\ell)\|_p^{p}\qquad\text{and}\qquad\sum_{j\ge1}|\lambda_j(T)|^{p}=\sum_\ell\sum_{j\ge1}|\lambda_j(T_\ell)|^{p}.
		\end{equation*}
		By Theorem \ref{thm:products}, applied to $T_\ell$, the sequence of singular values of $\Delta_t^{n+1}(T_\ell)$ is weakly log-majorized by that of $\Delta_t^{n}(T_\ell)$ for every $n$, and \cite[Corollary 1.10]{Simon}, applied with the function $u\mapsto u^{p}$, turns this into $\|\Delta_t^{n+1}(T_\ell)\|_p^{p}\le\|\Delta_t^{n}(T_\ell)\|_p^{p}$. Thus each sequence $\|\Delta_t^{n}(T_\ell)\|_p^{p}$ is non-increasing in $n$, so it is dominated by its first term $\|T_\ell\|_p^{p}$, and $\sum_\ell\|T_\ell\|_p^{p}=\|T\|_p^{p}<\infty$. Suppose that
		\begin{equation}\label{eq:summand}
			\lim_{n\to\infty}\|\Delta_t^{n}(T_\ell)\|_p^{p}=\sum_{j\ge1}|\lambda_j(T_\ell)|^{p}
		\end{equation}
		holds for every $\ell$. Then dominated convergence, applied to the sum over $\ell$, gives
		\begin{align*}
			\lim_{n\to\infty}\|\Delta_t^{n}(T)\|_p^{p}&=\lim_{n\to\infty}\sum_\ell\|\Delta_t^{n}(T_\ell)\|_p^{p}=\sum_\ell\lim_{n\to\infty}\|\Delta_t^{n}(T_\ell)\|_p^{p}\\&=\sum_\ell\sum_{j\ge1}|\lambda_j(T_\ell)|^{p}=\sum_{j\ge1}|\lambda_j(T)|^{p},
		\end{align*}
		and taking $p$-th roots gives the assertion. It therefore suffices to prove \eqref{eq:summand} for every $\ell$.
		
		If $T_\ell$ acts on a finite-dimensional space, then \eqref{eq:summand} is a finite sum of the limits in Theorem \ref{thm:main-intro}. If $T_\ell$ is a unilateral weighted shift, then $T_\ell\in\mathcal{S}_p$, since $\|T_\ell\|_p\le\|T\|_p$, so the left-hand side of \eqref{eq:summand} is $0$ by Proposition \ref{prop:shift}, and the right-hand side is $0$ because $T_\ell$ has no non-zero eigenvalues.
	\end{proof}
	
	\bigskip
	\section{The Duggal transform}\label{sec:duggal}
	\bigskip
	
	We conclude with the endpoint $t=1$, which is excluded from Theorem \ref{thm:main-intro}. The Duggal transform is of a different nature than the transforms $\Delta_t$ with $0<t<1$: up to unitary equivalence, it only cuts an operator down to the closure of its range. For $T\in\mathfrak{B}(\mathcal{H})$ and $n\ge0$ we put $\mathcal{R}_n=\overline{\operatorname{ran}T^{n}}$, so that $\mathcal{H}=\mathcal{R}_0\supseteq\mathcal{R}_1\supseteq\mathcal{R}_2\supseteq\cdots$, and we let $\mathcal{R}_\infty=\bigcap_{n\ge0}\mathcal{R}_n$. Each $\mathcal{R}_n$ is invariant under $T$, since $T\mathcal{R}_n\subseteq\overline{T\operatorname{ran}T^{n}}=\mathcal{R}_{n+1}$ by the continuity of $T$, and so is $\mathcal{R}_\infty$. We write $Q_n$ and $Q_\infty$ for the orthogonal projections of $\mathcal{H}$ onto $\mathcal{R}_n$ and $\mathcal{R}_\infty$.
	
	\begin{proposition}\label{prop:duggal}
		Let $T\in\mathfrak{B}(\mathcal{H})$ and $n\in\mathbb{N}$. Then there is a Hilbert space $\mathcal{Z}$, possibly $\{0\}$, such that $\Delta_1^{n}(T)$ is unitarily equivalent to the operator $T|_{\mathcal{R}_n}\oplus0$ on $\mathcal{R}_n\oplus\mathcal{Z}$. If $T$ is compact, then $s_j(\Delta_1^{n}(T))=s_j(TQ_n)$ for every $j$, and
		\begin{equation*}
			\lim_{n\to\infty}s_j(\Delta_1^{n}(T))=s_j(TQ_\infty)\qquad\text{for every }j.
		\end{equation*}
	\end{proposition}
	
	\begin{proof}
		Let $T=U|T|$ be the polar decomposition. Since $U^*U$ is the orthogonal projection onto $(\ker T)^\perp=\overline{\operatorname{ran}|T|}$, we have $U^*U|T|=|T|$, and therefore
		\begin{equation*}
			\Delta_1(T)=|T|U=U^*U|T|U=U^*TU,
		\end{equation*}
		an identity contained in \cite[Lemma 2.3]{FJKP}. The partial isometry $U$ vanishes on $\ker T$ and maps $(\ker T)^\perp$ isometrically onto its final space $\overline{\operatorname{ran}T}=\mathcal{R}_1$. Let $W\colon(\ker T)^\perp\to\mathcal{R}_1$ be the resulting unitary operator. For $y=Wx$ with $x\in(\ker T)^\perp$ we have $U^*y=U^*Ux=x$, so $U^*$ agrees with $W^{-1}$ on $\mathcal{R}_1$. Now let $x\in\mathcal{H}$ and write $x=x_1+x_2$ with $x_1\in(\ker T)^\perp$ and $x_2\in\ker T$. Then $Ux=Wx_1\in\mathcal{R}_1$, hence $TUx=TWx_1\in\mathcal{R}_1$ because $\mathcal{R}_1$ is invariant under $T$, and consequently
		\begin{equation*}
			\Delta_1(T)x=U^*TUx=W^{-1}\big(T|_{\mathcal{R}_1}\big)Wx_1\in(\ker T)^\perp.
		\end{equation*}
		Thus $\Delta_1(T)$ vanishes on $\ker T$, leaves $(\ker T)^\perp$ invariant, and acts there as $W^{-1}(T|_{\mathcal{R}_1})W$. The operator $V=W\oplus I_{\ker T}$ is unitary from $\mathcal{H}=(\ker T)^\perp\oplus\ker T$ onto $\mathcal{R}_1\oplus\ker T$, and $V\Delta_1(T)V^{-1}=T|_{\mathcal{R}_1}\oplus0$. This is the first assertion for $n=1$, with $\mathcal{Z}=\ker T$. The argument applies to every bounded operator on every Hilbert space.
		
		Suppose that the first assertion holds for some $n$, say $V\Delta_1^{n}(T)V^{-1}=T|_{\mathcal{R}_n}\oplus0$ with a unitary operator $V$ from $\mathcal{H}$ onto $\mathcal{R}_n\oplus\mathcal{Z}$. The polar decomposition is compatible with unitary equivalence and with orthogonal direct sums, hence so is $\Delta_1$, and $\Delta_1(0)=0$. Therefore
		\begin{equation*}
			V\Delta_1^{n+1}(T)V^{-1}=\Delta_1\big(T|_{\mathcal{R}_n}\oplus0\big)=\Delta_1\big(T|_{\mathcal{R}_n}\big)\oplus0.
		\end{equation*}
		By the case $n=1$, applied to the operator $T|_{\mathcal{R}_n}$ on the Hilbert space $\mathcal{R}_n$, the operator $\Delta_1(T|_{\mathcal{R}_n})$ is unitarily equivalent to the restriction of $T|_{\mathcal{R}_n}$ to the closure of its range, padded by a zero operator. The range of $T|_{\mathcal{R}_n}$ is $T\mathcal{R}_n$, and
		\begin{equation*}
			\operatorname{ran}T^{n+1}=T\operatorname{ran}T^{n}\subseteq T\mathcal{R}_n\subseteq\overline{T\operatorname{ran}T^{n}}=\mathcal{R}_{n+1},
		\end{equation*}
		so that $\overline{T\mathcal{R}_n}=\mathcal{R}_{n+1}$ and the restriction in question is $T|_{\mathcal{R}_{n+1}}$. Hence $\Delta_1^{n+1}(T)$ is unitarily equivalent to $T|_{\mathcal{R}_{n+1}}\oplus0$, where the zero operator now acts on the orthogonal direct sum of $\ker(T|_{\mathcal{R}_n})$ and $\mathcal{Z}$. This completes the induction.
		
		Let now $T$ be compact. Then $\Delta_1^{n}(T)$ is compact, being unitarily equivalent to $T|_{\mathcal{R}_n}\oplus0$. With respect to $\mathcal{H}=\mathcal{R}_n\oplus\mathcal{R}_n^{\perp}$ we have $TQ_n=T|_{\mathcal{R}_n}\oplus0$, again because $\mathcal{R}_n$ is invariant under $T$. Since $|A\oplus0|=|A|\oplus0$ and $|VAV^{-1}|=V|A|V^{-1}$ for every unitary operator $V$, both $|TQ_n|$ and $|\Delta_1^{n}(T)|$ are unitarily equivalent to the orthogonal direct sum of $|T|_{\mathcal{R}_n}|$ and a zero operator. They have therefore the same non-zero eigenvalues with the same multiplicities, and by the convention of Section \ref{sec:prelim} this means that $s_j(\Delta_1^{n}(T))=s_j(TQ_n)$ for every $j$.
		
		It remains to prove the limit relation, and by what has just been shown it suffices to prove that $s_j(TQ_n)\to s_j(TQ_\infty)$ for every $j$. We first show that $Q_n\to Q_\infty$ strongly. Let $x\in\mathcal{H}$. For $m>n$ the operator $Q_n-Q_m$ is the orthogonal projection onto $\mathcal{R}_n\ominus\mathcal{R}_m$, so that
		\begin{equation*}
			\|Q_nx-Q_mx\|^2=\|Q_nx\|^2-\|Q_mx\|^2.
		\end{equation*}
		The sequence $(\|Q_nx\|^2)_{n\ge0}$ is non-increasing and bounded below, hence convergent, so $(Q_nx)_{n\ge0}$ is a Cauchy sequence. Let $y$ be its limit. For every $m$ the vectors $Q_nx$ with $n\ge m$ lie in the closed subspace $\mathcal{R}_m$, hence so does $y$. Thus $y\in\mathcal{R}_\infty$ and $Q_\infty y=y$. On the other hand $Q_\infty Q_n=Q_\infty$ for every $n$, so that $Q_\infty y=\lim_nQ_\infty Q_nx=Q_\infty x$. Hence $Q_nx\to Q_\infty x$.
		
		Next, $\|TQ_n-TQ_\infty\|\to0$. Otherwise there are $\delta>0$, a strictly increasing sequence $(n_i)_{i\ge1}$ of positive integers and unit vectors $x_i\in\mathcal{H}$ with $\|T(Q_{n_i}-Q_\infty)x_i\|\ge\delta$ for every $i$. The vectors $z_i=(Q_{n_i}-Q_\infty)x_i$ are bounded, and for every $y\in\mathcal{H}$ we have
		\begin{equation*}
			|\langle z_i,y\rangle|=|\langle x_i,(Q_{n_i}-Q_\infty)y\rangle|\le\|(Q_{n_i}-Q_\infty)y\|\longrightarrow0,
		\end{equation*}
		so $z_i\to0$ weakly. Since $T$ is compact, $Tz_i\to0$ in norm \cite[Theorem VI.11]{RS}, which contradicts $\|Tz_i\|\ge\delta$.
		
		Let $k\in\mathbb{N}$ with $k\le d$. The map $S\mapsto\wedge^kS$ is norm continuous on $\mathfrak{B}(\mathcal{H})$, so Lemma \ref{lem:wedge-compact} gives
		\begin{equation}\label{eq:duggal-products}
			\prod_{j=1}^{k}s_j(TQ_n)=\|\wedge^k(TQ_n)\|\longrightarrow\|\wedge^k(TQ_\infty)\|=\prod_{j=1}^{k}s_j(TQ_\infty)\qquad(n\to\infty),
		\end{equation}
		and for $k=0$ we interpret both products in \eqref{eq:duggal-products} as $1$. Put $\rho=\rho(|TQ_\infty|)$, the number of non-zero singular values of $TQ_\infty$, and let $j$ be a positive integer. If $j\le\rho$, then $\prod_{i<j}s_i(TQ_\infty)>0$, so by \eqref{eq:duggal-products} the products $\prod_{i<j}s_i(TQ_n)$ are positive for all large $n$, and dividing \eqref{eq:duggal-products} for $k=j$ by \eqref{eq:duggal-products} for $k=j-1$ gives $s_j(TQ_n)\to s_j(TQ_\infty)$. If $j>\rho$, then $\rho<j\le d$, so \eqref{eq:duggal-products} is available for $k=\rho+1$ and $k=\rho$. Since $\prod_{i\le\rho+1}s_i(TQ_\infty)=0$ while $\prod_{i\le\rho}s_i(TQ_\infty)>0$, the same division gives $s_{\rho+1}(TQ_n)\to0$, and as $s_j(TQ_n)\le s_{\rho+1}(TQ_n)$ we obtain $s_j(TQ_n)\to0=s_j(TQ_\infty)$.
	\end{proof}
	
	\begin{remark}
		Proposition \ref{prop:duggal} shows why the endpoint $t=1$ is excluded from Theorem \ref{thm:main-intro}. If $T$ has dense range, then $\mathcal{R}_n=\mathcal{H}$ for every $n$, so that all Duggal iterates of $T$ are unitarily equivalent to $T\oplus0$, and if $T$ is moreover injective, they are unitarily equivalent to $T$. This applies to every invertible matrix. For
		\begin{equation*}
			T=\begin{pmatrix}1&1\\0&1\end{pmatrix}
		\end{equation*}
		one obtains $s_1(\Delta_1^{n}(T))=(1+\sqrt5)/2$ for every $n$, whereas $|\lambda_1(T)|=1$, and the same happens for the rank two operator $T\oplus0$ on any Hilbert space. Thus the conclusion of Theorem \ref{thm:main-intro} fails for $t=1$ already for operators of finite rank, and so does \eqref{eq:yamazaki}, in accordance with the example in \cite{BCLY}. In general, the singular values of the Duggal iterates of a compact operator converge to the singular values of $TQ_\infty$ rather than to the moduli of the eigenvalues of $T$, and the two coincide for normal $T$, since then $\mathcal{R}_n=(\ker T)^\perp$ for every $n\ge1$ and $TQ_\infty=T$.
	\end{remark}
	
	\begin{remark}\label{rem:separable}
		The separability of $\mathcal{H}$, assumed throughout, is a matter of convenience, since it is the setting of the sources \cite{GK} and \cite{Simon} on which we rely. All results of this paper remain valid for compact operators on an arbitrary complex Hilbert space $\mathcal{H}$. Indeed, let $T\in\mathcal{K}(\mathcal{H})$. For $T=0$ all assertions are trivial, so let $T\neq0$. The image of the unit ball under $T$ is relatively compact, hence separable, so $\operatorname{ran}T$ is separable, and the same holds for $\operatorname{ran}T^*$ because $T^*$ is compact. Thus $\mathcal{H}_0=\overline{\operatorname{ran}T+\operatorname{ran}T^*}$ is a separable closed subspace, it is invariant under $T$ and $T^*$, and its orthogonal complement is $\ker T^*\cap\ker T$. Hence $T=T_0\oplus0$ with respect to $\mathcal{H}=\mathcal{H}_0\oplus\mathcal{H}_0^{\perp}$, where $T_0=T|_{\mathcal{H}_0}\in\mathcal{K}(\mathcal{H}_0)$. The polar decomposition and the functional calculus are compatible with orthogonal direct sums, so $\Delta_t^{n}(T)=\Delta_t^{n}(T_0)\oplus0$ for all $t\in[0,1]$ and $n\ge0$, and $\overline{\operatorname{ran}T^{n}}=\overline{\operatorname{ran}T_0^{n}}$ for every $n\ge1$. The non-zero singular values and the non-zero eigenvalues of an operator of the form $A\oplus0$, counted with multiplicity, are those of $A$. Consequently every quantity appearing in the results above, computed for $T$, coincides with the same quantity computed for $T_0$ on the separable space $\mathcal{H}_0$, the sequences of singular values and eigenvalues being extended by zeros when $\mathcal{H}_0$ is finite-dimensional, and the results for $T$ follow from those for $T_0$.
	\end{remark}

	\bigskip
	\section*{Declarations}
	\noindent{\bf{Funding}}\\
	This work has been supported by the Ministry of Science, Technological Development and Innovation of the Republic of Serbia [Grant Number: 451-03-34/2026-03/200102].
	
	\vspace{0.5cm}
	
	\noindent{\bf{Availability of data and materials}}\\
	\noindent No data were used to support this study.
	\vspace{0.5cm}\\
	\noindent{\bf{Competing interests}}\\
	\noindent The author declares that he has no competing interests.
	\vspace{0.5cm}\\
	\noindent{\bf{Declaration of generative AI and AI-assisted technologies in the manuscript preparation process}}\\
	\noindent During the preparation of this work, the author used Claude (Anthropic) to identify relevant references, to verify the mathematical arguments, to assist in the derivation and drafting of parts of the results, and for language editing and improvement of the exposition. The proof of the second formula in Lemma \ref{lem:wedge-compact}, for which no direct proof could be located in the literature, was obtained with the help of Claude (Anthropic). ChatGPT (OpenAI) was additionally consulted for pre-submission checks. The author reviewed, verified and edited all generated content and takes full responsibility for the content of the published article.


\end{document}